\documentclass{article}
\usepackage[utf8]{inputenc}
\usepackage[english]{babel}
\usepackage{amsmath,amssymb,euscript,amsthm,amsfonts,mathrsfs,amscd,mathtext}
\usepackage{color}
\usepackage{cite,enumerate,float,indentfirst}
\usepackage{graphicx}
\usepackage{setspace}
\usepackage{mathrsfs} 

\usepackage[left=20mm, top=20mm, right=20mm, bottom=20mm, footskip=10mm, nohead, nomarginpar, driver=xetex]{geometry}

\theoremstyle{plain}
\newtheorem{thm}{Theorem}

\theoremstyle{definition}

\theoremstyle{remark}
\newtheorem{rem}{Remark}

\newcommand{\bb}{\color{blue}}
\newcommand{\gr}{\color{green}}
\newcommand{\ww}{\color{white}}
\newcommand{\rr}{\color{red}}

\newcommand{\mg}{\color{magenta}}
\newcommand{\dd}{\ensuremath{\displaystyle}}

\newcommand{\supl}{\ensuremath{\sup\limits}}
\newcommand\intl{\ensuremath{\int\limits}}
\newcommand{\suml}{\ensuremath{\sum\limits}}
\newcommand{\PP}{\ensuremath{\mathbf{P}}}
\newcommand{\UUU}{\ensuremath{\mathscr{U}}}

\newcommand{\PPP}{\ensuremath{\mathscr{P}}}
\newcommand{\XXX}{\ensuremath{\mathscr{X}}}
\newcommand{\BBB}{\ensuremath{\mathscr{B}}}
\newcommand{\SSS}{\ensuremath{\mathscr{S}}}

\newcommand{\EE}{\mathbb E}

\newcommand{\ud}{\,\mathrm{d}\,}
\newcommand*{\TR}{\hfill\ensuremath{\triangleright}}
\newcommand{\bd}{\stackrel{ {\rm def}}{=\!\!\!=}}
\title{On some quasi-regenerative reliability system}
\author{Galina Zverkina\footnote{V. A. Trapeznikov Institute of Control Sciences of Russian Academy of Sciences}~\footnote{The work is supported by RFBR grant No 20-01-00575A}}
\date{December 2022}

\begin{document}

\maketitle

\section{Introduction}

We consider the reliability system consisting of two {\it dependent} restorable elements. 
Both restorable elements have the same characteristics.
However, the ``speed'' of work or repair depends on the {\it full state} of the system.

At first, the first (main) element works, and the serviceable second element does not work at full capacity - it is in a warm reserve.

If the main element fails, then the reserve element starts working at full speed -- it is now the main one.

Or if a working element has already been working for a long time, then the second element begins to warm up in order to start working faster when the first element fails -- to immediately switch to the active mode of operation when a working element fails.

In addition, if the element works for a very long time, its reliability decreases and the likelihood of its failure in the near future increases.

In addition, if both elements are faulty, then they are repairing at a faster rate than in another situation.

The probability of failure of a working element or the probability of repair of a repaired element is determined by {\it intensities} that depend only on the {\it full state} of the system under study.

The full state of this reliability system at the time $t$ is a vector  $X_t\bd(n_1,x_1;n_2,x_2)\big[=(n_1(t),x_1(t);n_2(t),x_2(t))\big]$, where:
$$n_j(t)\bd\left\{
\begin{array}{ll}
1,&\mbox{if $j$-th element is in working state at the time $t$};
\\ 
0,&\mbox{if $j$-th element is in failure state at the time $t$};
\end{array}
\right.
$$
the variable $x_j(t)$ is the elapsed time of the stay of $j$-th element in the status $n_j  $ at the time $t$.

The behaviour of $X_t$, its distribution give the full information about this reliability system.

The probability of a status $n_j(t)$ change over a time interval $(t,t+\Delta)$ is equal to  $\lambda_j(X_t)\Delta+o(\Delta)$; $\lambda_j(X_t)$ is an intensity of  a status $n_j(t)$ change at the time $t$.

Under these assumptions, the process $X_t$ is Markov on the state space $\XXX\bd\{(0;1)\times \mathbb{R}_+\times(0;1)\times \mathbb{R}_+\}$.

\begin{figure}[h]
    \centering
\begin{picture}(400,120)
\put(-10,35){\small $\hat t_0=0$}
\put(0,30){\vector(1,0){400}}
\put(0,100){\vector(1,0){400}}
\put(0,30){\circle*{3}}
\put(0,100){\circle*{3}}
\put(0,15){ELEMENT II.}
\thicklines
\qbezier(0,30)(50,40)(100,30)
\put(30,40){\small work [1]}
\put(100,35){\small $\hat t_1$}
\qbezier(100,30)(120,20)(140,30)
\put(100,15){\small repair [0]}
\put(140,35){\small $\hat t_2$}
\qbezier(140,30)(170,40)(200,30)
\put(147,40){\small work [1]}
\put(200,35){\small $\hat t_3$}
\qbezier(200,30)(230,20)(260,30)
\put(210,15){\small repair [0]}
\put(260,35){\small $\hat t_4$}
\qbezier(260,30)(297,40)(330,30)
\put(282,40){\small work [1]}
\put(320,35){\small $\hat t_5$}
\qbezier(330,30)(345,20)(360,30)
\put(345,15){\small repair [0]}
\put(355,35){\small $\hat t_6$}
\multiput(365,35)(10,0){5}{\circle*{2}}
\put(-10,105){\small $t_0=0$}
\put(0,85){ELEMENT I.}
\qbezier(0,100)(80,110)(160,100)
\put(30,110){\small work [1]}
\put(150,105){\small $ t_1$}
\qbezier(160,100)(200,90)(240,100)
\put(185,85){\small repair [0]}
\put(240,105){\small $ t_2$}
\qbezier(240,100)(280,110)(320,100)
\put(245,110){\small work [1]}
\put(320,105){\small $ t_3$}
\qbezier(320,100)(350,90)(380,100)
\put(343,85){\small repair [0]}
\put(375,105){\small $ t_4$}
\multiput(385,105)(10,0){3}{\circle*{2}}
\multiput(180, 15)(0,5){21}{\bb \line(0,1){3}}
\put(182,20){\small $\theta_1$}
\put(182,105){\small $\theta_1$}
\put(160,100){\gr\line (0,1){15}}
\put(180,110){\gr \vector(-1,0){20}}
\put(168,113){\small $x_1$}
\put(140,30){\gr\line (0,-1){15}}
\put(180,15){\gr \vector(-1,0){40}}
\put(160,20){\small $y_1$}
\multiput(280, 15)(0,5){21}{\bb \line(0,1){3}}
\put(240,100){\gr\line (0,-1){15}}
\put(282,20){\small $\theta_2$}
\put(282,92){\small $\theta_2$}
\put(280,85){\gr \vector(-1,0){40}}
\put(260,90){\small $x_2$}
\put(260,30){\gr\line (0,-1){15}}
\put(280,15){\gr \vector(-1,0){20}}
\put(265,20){\small $y_2$}
\multiput(340, 15)(0,5){21}{\bb \line(0,1){3}}
\put(342,32){\small $\theta_3$}
\put(342,102){\small $\theta_3$}
\put(320,100){\gr\line (0,1){15}}
\put(340,115){\gr \vector(-1,0){20}}
\put(325,120){\small $x_3$}
\put(330,30){\gr\line (0,1){15}}
\put(340,45){\gr \vector(-1,0){10}}
\put(330,50){\small $y_3$}
\end{picture}
\caption{This scheme is a visualization of the work and repair of a reliability system with a warm reserve. Here at the time $\theta_1$ we have $X_{\theta_1}=(0,x_1;1,y_1)$, and the main element is II; at the time $\theta_2$ we have $X_{\theta_2}=(1,x_2;1,y_2)$, and the main element is I; at the time $\theta_2$ we have  $X_{\theta_3}=(0,x_3;0,y_3)$, and the system in failure state.}
\end{figure}
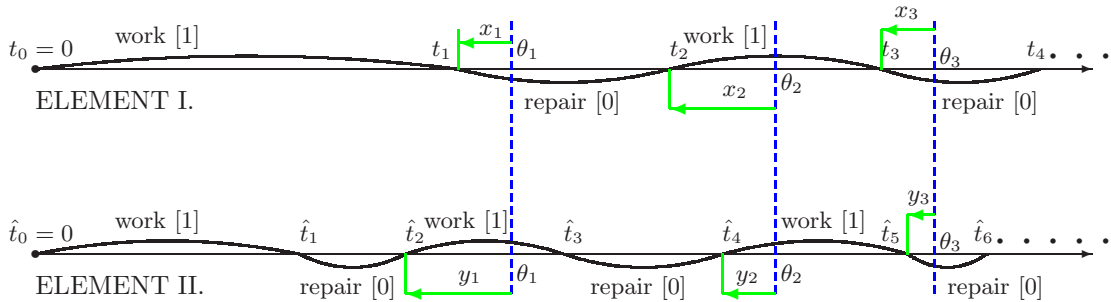

In our notation, the set $\SSS_{0,0}\bd \{(0,x;0,y), \; x,y\in \mathbb{R}_{+}\}$ corresponds to the failure state of the system. 
And the set $\SSS_{>0}\bd \{(n_1,x_1;n_2,x_2): \; n_1+n_2>0\}$ corresponds to the system operating state, and $\PP\{X_t\in \SSS_{>0}\}$ is equal to the availability factor of considered system.

Therefore, knowledge of the probability distribution of the stochastic process $X_t$ is important for calculating the characteristics of such systems.

\subsection{Reliability system considered earlier}

We suppose that the intensities of failure for working elements and the intensities of repair for repaired elements depend on the process $X_t$: $\lambda_i^{(1)}=\lambda_i^{(1)}(X_t)$ is the intensities of failure of $i$-th element, and $\lambda_i^{(0)}=\lambda_i^{(0)}(X_t)$ is the intensities of repair of $i$-th element.
In this situation, the process $X_t$ is Markov on the state space $\XXX\bd(\{0;1\}\times \mathbb R_{+})^2$ with a standard $\sigma$-algebra.
If 
\begin{equation}\label{exp}
   \lambda_i^{(n)}(X_t)\equiv \mathrm{const}, 
\end{equation} then this is a classic case with exponential d.f.'s.

If $0<c\leqslant\lambda_i^{(n)}(X_t)\leqslant C<\infty$, then this ``quasi-exponential'' case, it is studied in the paper (using some properties of work and repair times following from (\ref{exp}) -- see \cite{GZ}).

In the paper \cite{ver_av}, a reliability system with a warm reserve was studied under the following conditions for intensities. 
There exist constants $\gamma, \Gamma>0$ such that for any $Z =(i,x;j,y)\in \XXX$ (state space)
\begin{equation}\label{usl}
\begin{array}{c}
0<\dd\frac{\gamma}{1+x} \leq \lambda(Z) \le \Gamma<\infty;
\\
\\
0<\dd\frac{\gamma}{1+y}\leq {\mu}(Z)\le \Gamma < \infty.
\end{array}
\end{equation}
If the conditions (\ref{usl}) are satisfied, then the process $X_t$ is recurrent with finite return time to some compact. 
Therefore, the process $X_t$ is ergodic.

In both cases (\cite{GZ,ver_av}), by default, it is assumed that the switching time (work-repair and repair-work) is instantaneous, and the distribution of work and repair time is absolutely continuous.

\section{Our goal}
Thus,  we will study the situation when switching between operation and repair modes is not instantaneous, and the distribution of work and repair time can have a discrete part.
In addition to this, we will study the case when all switching times are random variables bounded above by a constant with a probability of 1.


\subsection{Intensities}

First, we define the generalized intensity for random variables whose distribution has a discrete component.

In order to explore the situation of mixed random variables (i.e. the situation where the distribution can have discrete and continuous components), we will use the concept of generalized intensity introduced in the paper \cite{KZ}.

Recall, that the intensity of the distribution of a continuous positive random variable (r.v.) $\xi$ with distribution function (d.f.) $F(s)$ is the function
$
\dd\lambda(s)\bd \frac{F'(s)}{1-F(s)},
$
and $\PP\{\xi\in(s,s+\Delta )|\xi>s\}=\lambda(s)\Delta+o(\Delta)$.

The function $\lambda(s)$ defines d.f. $F(s)$:
\begin{equation} \label{defint}
   F(s)=1-\dd e^{-\intl_0^s \lambda(v) \ud v}.    
\end{equation}

If the distribution of $\xi$ is mixed (non-singular!!!), i.e. d.f. $F(s)$ has jumps, then we can put
\begin{equation}\label{genint}
    \begin{array}{l}
  f(s)=\left
\{
\begin{array}{ll}
F'(s), & \mbox{if } F'(s) \mbox{ exists; }
 \\ \\
0, & \mbox{otherwise};
\end{array}
\right.     \\ \\
\lambda(s)\bd \dd\frac{f(s)}{1-F(s)}-\suml_{i}\delta(s-a_i)\ln(F(a_i+0)-F(a_i-0)),
\end{array}
\end{equation}
where $\{a_i\}$ is the set of discontinuity points of $F(s)$, and $\delta(s)$ is a standard $\delta$-function;
\\
The formula (1) is true with these notations.

{\it 
So, the distribution of  non-negative r.v.'s -- continuous, discrete and mixed -- can be defined by intensities.}

However, the absolutely arbitrary dependence of the intensity of the change in the base state (work, repair) does not allow us to study the behavior of the reliability system with a warm reserve.

\subsection{Condition for intensities}

In the future, for the time of operation and / or repair, we will use the conditions {\bf a}--{\bf d} or {\bf a}--{\bf d}.

{\it The following conditions for the (generalized) intensities $\lambda_i^{(n)}(s)$ are assumed:}
\begin{enumerate}
\item [\bf a.]{\it The (generalized) measurable non-negative functions $ \varphi(s)$ and $Q(s)$ 
exist 
such that for all $s\geqslant 0$, $${ \varphi(s)\leqslant\lambda^{(n)}_i(s) \leqslant Q(s)};$$}
\item [\bf b.]{\it $\dd\intl_0^\infty  \varphi(s) \ud s = \infty$, and  $\dd\intl_0^\infty x^{k-1}  \exp\left(-\intl_0^x  \varphi(s)\ud s\right)\ud x<\infty $ for some $k\geqslant 2$;}
\item [\bf c.]{\it There is a neighborhood of zero $\UUU\bd(-\varepsilon,\varepsilon)$ such that $\dd\intl_\UUU Q(s)\ud s<1$;
\item [\bf d.]There exists the constant $T\geqslant 0$ such that $ \varphi(s)>0$ a.s. for all $s>T$.}
\end{enumerate}
\begin{rem}
  Condition {\bf (a)} holds: $G(s)=\PP\{\zeta\leqslant s\}\geqslant F_i^{(n)}(s)= \PP\{\xi^{(n)}_i\leqslant s\} \geqslant \Phi(s) =  \PP\{\widetilde\zeta\leqslant s\}$, or $\zeta\succsim \xi^{(n)}_i \succsim \widetilde\zeta$ are ordered in distribution $\big[$here $G(s)$ has intensity $Q(s)$, $\Phi(s)$ has intensity $\varphi(s)\big]$.
\TR  
\end{rem}

\begin{rem}
 Condition {\bf (b)} holds: there exists $\EE\,\zeta^k<\infty \Rightarrow  \EE\,\widetilde\zeta^k<\infty$ and $\EE\,\left(\xi^{(n)}_i\right)^k<\infty$.
 \TR   
\end{rem}

\begin{rem}
    Condition {\bf (c)} holds: $\EE\,\zeta^2>0$, and $\mathrm{Var}(\zeta)>0$.
 \TR
\end{rem} 

\begin{rem}
Condition { \bf (d)} holds: $\Phi'(x)>0$ a.s. for $x>T$, i.e. we may consider delayed switching for the considered reliability system.
 \TR
\end{rem}

\subsection{The studied warm standby systems is non-regenerative}

An analysis of the behaviour of both elements of the investigated reliability system shows that the periods of operation and repair of both elements depend on each other.
Consequently, for each individual element of this reliability system there is no independence between the periods of work and repair.

\begin{figure}[h]
    \centering
\begin{picture}(400,80)
\put(-10,35){\small $\hat t_0=0$}
\put(0,30){\vector(1,0){400}}
\put(0,60){\vector(1,0){400}}
\put(0,30){\circle*{3}}
\put(0,60){\circle*{3}}
\put(0,20){ELEMENT II.}
\thicklines
\qbezier(0,30)(50,40)(100,30)
\put(30,40){\small work [1]}
\put(100,35){\small $\hat t_1$}
\qbezier(100,30)(120,20)(140,30)
\put(110,15){\small repair [0]}
\put(140,35){\small $\hat t_2$}
\qbezier(140,30)(170,40)(200,30)
\put(150,40){\small work [1]}
\put(200,35){\small $\hat t_3$}
\qbezier(200,30)(230,20)(260,30)
\put(210,15){\small repair [0]}
\put(260,35){\small $\hat t_4$}
\qbezier(260,30)(297,40)(330,30)
\put(270,40){\small work [1]}
\put(320,35){\small $\hat t_5$}
\qbezier(330,30)(345,20)(360,30)
\put(330,15){\small repair [0]}
\put(355,35){\small $\hat t_6$}
\multiput(365,35)(10,0){5}{\circle*{2}}
\put(-10,65){\small $ t_0=0$}
\put(0,50){ELEMENT I.}
\qbezier(0,60)(80,70)(160,60)
\put(30,70){\small work [1]}
\put(150,65){\small $ t_1$}
\qbezier(160,60)(200,50)(240,60)
\put(190,45){\small repair [0]}
\put(240,65){\small $ t_2$}
\qbezier(240,60)(280,70)(320,60)
\put(250,70){\small work [1]}
\put(320,65){\small $ t_3$}
\qbezier(320,60)(350,50)(380,60)
\put(330,45){\small repair [0]}
\put(375,65){\small $ t_4$}
\multiput(385,65)(10,0){3}{\circle*{2}}
\put(30,35){\bb\vector(1,1){30}}
\put(70,65){\bb\vector(-1,-1){30}}
\put(110,27){\bb\vector(1,2){18}}
\put(135,62){\bb\vector(-1,-2){18}}
\put(155,35){\bb\vector(1,2){12}}
\put(175,58){\bb\vector(-1,-2){12}}
\end{picture}
\caption{Dependence between the work and repair periods of warm standby systems.}
\end{figure}
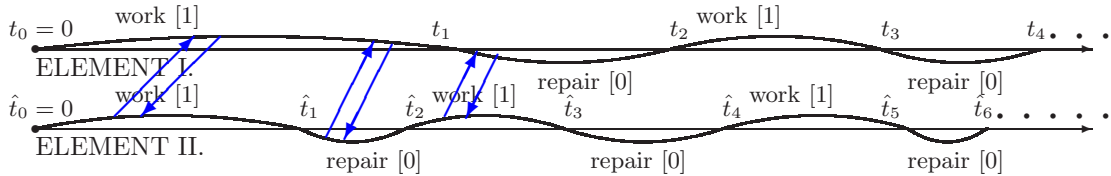

The sequence of periods $(t_0;t_1)$, $(t_1;t_2)$, $(t_2;t_3)$, $(t_3;t_4)$, \ldots and  $(\hat t_0;\hat t_1)$, $(\hat t_1;\hat t_2)$, $(\hat t_2;\hat t_3)$, $(\hat t_3;\hat t_4)$, \ldots ``looks like'' an alternating renewal process.

Through the dependence of the periods $(t_0;t_1)$ and $(\hat t_0;\hat t_1)$; $(t_0;t_1)$ and $(\hat t_1;\hat t_2)$; $(t_1;t_2)$ and $(\hat t_2;\hat t_3)$; etc., we have some dependencies between the lengths of the periods $(t_0;t_1)$, $(t_1;t_2)$, $(t_2;t_3)$, $(t_3;t_4)$, \ldots and between the lengths of the periods $(\hat t_0;\hat t_1)$, $(\hat t_1;\hat t_2)$, $(\hat t_2;\hat t_3)$, $(\hat t_3;\hat t_4)$, \ldots

For study, these dependencies cannot be arbitrary; below we indicate the conditions for the intensities.

These dependencies are ``weak'' in some sense.
Thus, the behaviour of each of the elements is described by an alternating quasi-renewal process, which is ``close'' to the classical renewal alternating process.

\subsection{How is the behaviour of this non-regenerative system studied.}

Thus, the first step of analysis of the behaviour of our system is the proof of the fact:
\begin{thm}
The process $X_t$ is quasi-regenerative, i.e. there exists regenerative Markov process $\tilde X_t$ (with finite regeneration time) on the state space $\XXX$ such that at any time the marginal distribution of $X_t$ is equal the marginal distribution of $\tilde X_t$.

Thus, the process $X_t$ is ergodic. \TR
\end{thm}

The proof based on the coupling method$^1$. 

The next step is the analysis of the behaviour of two processes $X_t$ and $X_t'$ with different initial state, but with the same transition probabilities.

By construction of {\it successful coupling} given by D. Griffeath \cite{zv5}, an upper bounds for the convergence rate of $X_t$ can be calculated.

\vspace{5mm}
Recall, that our goal is to compute an upper bound for  convergence rate of the distribution  of the reliability system in the case when the distributions of operating and recovery times may not be absolutely continuous, and mode switching may have a random {\it finite} delay.
The initial state of the system can be arbitrary.
 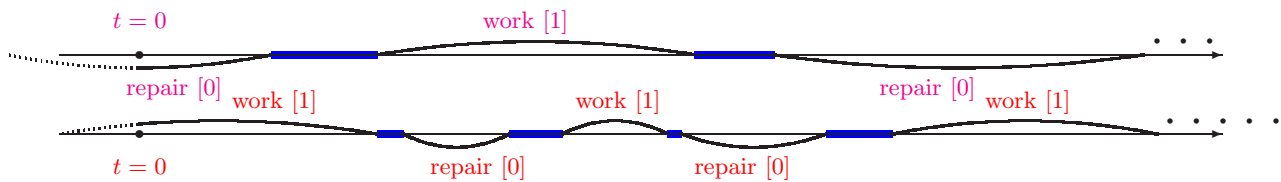
\begin{figure}[h]
    \centering
\begin{picture}(400,80)
\thicklines
\qbezier(-30,20)(30,30)(90,20)
\qbezier(-50,50)(0,40)(50,50)
\multiput(-50,0)(2,0){25}{\ww  \line(0,1){80}}
\put(35,30){\rr\small work [1]}
\put(90,20){\bb \line(1,0){10}}
\put(90,19){\bb \line(1,0){10}}
\put(90,21){\bb \line(1,0){10}}
\qbezier(100,20)(120,10)(140,20)
\put(110,5){\small\rr repair [0]}
\put(140,20){\bb \line(1,0){20}}
\put(140,19){\bb \line(1,0){20}}
\put(140,21){\bb \line(1,0){20}}
\qbezier(160,20)(180,30)(200,20)
\put(200,20){\bb \line(1,0){5}}
\put(200,19){\bb \line(1,0){5}}
\put(200,21){\bb \line(1,0){5}}
\put(165,30){\small \rr work [1]}
\qbezier(205,20)(232,10)(260,20)
\put(260,20){\bb \line(1,0){25}}
\put(260,19){\bb \line(1,0){25}}
\put(260,21){\bb \line(1,0){25}}
\put(210,5){\small \rr repair [0]}
\qbezier(285,20)(335,30)(385,20)
\put(320,30){\small\rr work [1]}
\multiput(390,25)(10,0){5}{\circle*{2}}
\put(-10,60){\small\mg $t=0$}
\put(-5,35){\small \mg repair [0]}
\put(50,50){\bb \line(1,0){40}}
\put(50,51){\bb \line(1,0){40}}
\put(50,49){\bb \line(1,0){40}}
\qbezier(90,50)(150,60)(210,50)
\put(130,60){\small\mg  work [1]}
\put(210,50){\bb \line(1,0){30}}
\put(210,51){\bb \line(1,0){30}}
\put(210,49){\bb \line(1,0){30}}
\qbezier(240,50)(310,40)(380,50)
\put(280,35){\small\mg repair [0]}
\multiput(385,55)(10,0){3}{\circle*{2}}
\put(-10,5){\small\rr $t=0$}
\thinlines
\put(-30,20){\vector(1,0){440}}
\put(-30,50){\vector(1,0){440}}
\put(0,20){\circle*{3}}
\put(0,50){\circle*{3}}
\end{picture}
\caption{Initial state of the process $X_t$ is $X_0=(0,x;1,y)$. Blue bars represent a switching time.}
\end{figure}

We suppose that the intensities of failure for working elements and the intensities of repair for repaired elements depend on the process $X_t$: $\lambda_i^{(1)}=\lambda_i^{(1)}(X_t)$ is the intensities of failure of $i$-th element, and $\lambda_i^{(0)}=\lambda_i^{(0)}(X_t)$ is the intensities of repair of $i$-th element.

\begin{rem}
For greater generality, we assume that all restrictions under the conditions of {\bf a}--{\bf d} are different for different elements.

So, we can give the following modified conditions
 
{\bf a'}. 
 {\it For all $s\geqslant 0$, $ \varphi_i^{(n)}(s)\leqslant\lambda_i^{(n)}(s) \leqslant Q_i^{(n)}(s)$;}
 
{\bf d'}. {\it There exists the constants $T_i^{(n)}\geqslant 0$ such that $ \varphi_i^{(n)}(s)>0$ a.s. for all $s>T_i^{(n)}$.}
  
  In this case, the upper bounds for $k$-th moments of $\xi^{(n)}_i$ are different, as well as its delay times.
  
  Put $C_j^{(n)}(\ell)\bd \dd\intl_0^\infty s^\ell\ud \Phi_j^{(n)}(s)$, where $\Phi_j^{(n)}(s)\bd 1-\dd e^{-\intl_0^s \varphi_j^{(n)}(v) \ud v}$ (see formula (\ref{defint})).
  
  The vector $\vec{C}(\ell)=\left\{C_1^{(1)}(\ell) ,C_2^{(1)}(\ell) ,C_1^{(0)}(\ell) ,C_2^{(0)}(\ell) \right\}$ depends on the vector  $\vec{\varphi}(s)=\left\{\varphi_1^{(1)}(s),\varphi_2^{(1)}(s),\varphi_1^{(0)}(s),\varphi_2^{(0)}(s)\right\}$.
  
  Also denote $\vec{Q}(s)=\left\{Q_1^{(1)}(s),Q_2^{(1)}(s),Q_1^{(0)}(s),Q_2^{(0)}(s)\right\}$.
\end{rem}

\section{Main results}

\begin{thm}\label{w-r}
Let the work and repair periods of both elements satisfy the conditions   {\bf a}--{\bf b}, and the work periods of both elements satisfy the  conditions  {\bf c}--{\bf d}.

Then the process $X_t$ is ergodic, i.e. the distribution $\PPP_t$ of $X_t$ at the time $t$, weak converges to the invariant stationary probabilistic distribution $\PPP$.\TR
\end{thm}

\begin{thm}\label{r-w}
Let the work and repair periods of both elements satisfy the conditions  {\bf a}--{\bf b}, and the repair periods of both elements satisfy the  conditions  {\bf c}--{\bf d}.

Then the process $X_t$ is ergodic, i.e. the distribution $\PPP_t$ of $X_t$ at the time $t$, weak converges to the invariant stationary probabilistic distribution $\PPP$.\TR
\end{thm}

\begin{thm}
In conditions of the Theorems \ref{w-r} and \ref{r-w}  for all $\ell\in(0,k-1]$ we give an algorithm of calculation of the constant $K(\ell,X_0,\vec{\varphi}(s),\vec{Q}(s))$ such that for all $t$,
 $$
 \|\PPP_t-\PPP\|_{TV}\leqslant \frac{K(\ell,X_0,\vec{\varphi}(s),\vec{Q}(s))}{t^\ell},$$
 where $\|\cdot - \cdot\|_{TV}$ is a total variation metric:
 $$
\|\PPP_t-\PPP\|_{TV}\bd \supl_{A\in \BBB(\XXX)}|\PPP_t(A)-\PPP(A)|.
$$  
\TR
\end{thm} 

\begin{thm}
Let the conditions of the Theorems \ref{w-r} and \ref{r-w} are satisfy, and let $\EE\exp(\alpha \tilde\zeta)<\infty$ for some $\alpha>0$.
In this case  we give an algorithm of calculation of the number $\beta\in(0,\alpha)$, and the constant $\tilde K(\beta,\ell,X_0,\vec{\varphi}(s),\vec{Q}(s))$ such that for all $t$,
 $$
 \|\PPP_t-\PPP\|_{TV}\leqslant\exp{-\beta t} \tilde K(\beta,\ell,X_0,\vec{\varphi}(s),\vec{Q}(s)),$$
 where $\|\cdot - \cdot\|_{TV}$ is a total variation metric:
 $$
\|\PPP_t-\PPP\|_{TV}\bd \supl_{A\in \BBB(\XXX)}|\PPP_t(A)-\PPP(A)|.
$$  
\TR
\end{thm} 

\subsection{1. Ergodicity}
The proof of ergodicity of the process $X_t$ based on the notion Generalized Markov Modulated Poisson Process (GMMPP) and generalized Lorden's inequality \cite{7,KZ}.

The consecutive pairs of periods of work and repair make up GMMPP that satisfy the conditions of ergodicity of such processes.

\subsection{2. Construction of upper bounds for convergence rate}

The process $X_t$ is quasi-regenerative, i.e. on some probability space there exists the Markov regenerative process $\widetilde{X}_t$ with the same {\it marginal} distributions as the process $X_t$: for all $A\in \BBB(\XXX)$ and for all $t\geqslant 0$, $\PP\{X_t\in A\}=\PP\{\widetilde{X}_t\in A\}$.

And for this Markov regenerative process $\widetilde{X}_t$ we can construct the constant $K(\ell,X_0,\vec{\varphi}(s),\vec{Q}(s))$.
The construction based on the {\it coupling method}. \cite{2,3}

These calculations are not optimal and can be improved by the use some properties of specific intensities.

\begin{rem}
The schema used in this work, can be used for reliability systems with many warm-reserve elements.

~

The construction of upper bounds for convergence rate is also possible for the case when at least one period (work or repair) for any element satisfies the condition {\bb \bf d}: for all $s>T$, $\varphi(s)>0$.

~

Naturally, an increase in the number of elements leads to an increase in the constant $K(\cdots)$.

~

Also it is important to understand that the coupling method does not provide an accurate estimate for stochastic processes in continuous time.
\end{rem}

\end{document}